\documentclass[12pt,reqno]{article}

\usepackage[usenames]{color}
\usepackage{amssymb}
\usepackage{amsmath}
\usepackage{amsthm}
\usepackage{amsfonts}
\usepackage{amscd}
\usepackage{graphicx}
\usepackage[colorlinks=true,
linkcolor=webgreen,
filecolor=webbrown,
citecolor=webgreen]{hyperref}

\definecolor{webgreen}{rgb}{0,.5,0}
\definecolor{webbrown}{rgb}{.6,0,0}

\usepackage{color}
\usepackage{fullpage}
\usepackage{float}
\usepackage{graphics}
\usepackage{latexsym}

\begin{document}

\theoremstyle{plain}
\newtheorem{theorem}{Theorem}
\newtheorem{corollary}[theorem]{Corollary}
\newtheorem{proposition}{Proposition}
\newtheorem{lemma}{Lemma}
\newtheorem{example}{Examples}
\newtheorem{remark}{Remark}

\begin{center}

{\large\bf
Forsyth--Glaisher-Type Series for $1/\pi$ Involving \\[6pt] Catalan Numbers and Harmonic Numbers}

\vskip 1.25cm

{\large
Kunle Adegoke \\
Department of Physics and Engineering Physics \\ Obafemi Awolowo University, 220005 Ile-Ife, Nigeria \\
\href{mailto:adegoke00@gmail.com}{\tt adegoke00@gmail.com}

\vskip 0.15 in

Robert Frontczak \\
Independent Researcher\\
 72764 Reutlingen, Germany \\
\href{mailto:robert.frontczak@web.de}{\tt robert.frontczak@web.de}

\vskip 0.15 in

Taras Goy \\
Faculty of Mathematics and Computer Science \\
Vasyl Stefanyk Carpathian National University \\
76018 Ivano-Frankivsk, Ukraine \\[4pt]
\href{mailto:taras.goy@pnu.edu.ua}{\tt taras.goy@pnu.edu.ua}
}

\end{center}

\vskip .15 in

\begin{abstract}
We develop an elementary method for constructing and evaluating families of Forsyth--Glaisher-type series for $1/\pi$. The approach is based on Newton's generalized binomial theorem, identities for generalized binomial coefficients, and Euler's beta integral, together with properties of the gamma and digamma functions. We establish several parameterized identities involving squares of central binomial coefficients and products of odd integers, from which we derive explicit evaluations in terms of $1/\pi$. By means of transformations connecting generalized binomial coefficients with Catalan numbers, we obtain further families of series involving squares of Catalan numbers. We also derive identities for series involving harmonic and odd harmonic numbers by differentiating suitable parameter-dependent formulas. These results yield evaluations containing both $1/\pi$ and logarithmic terms, particularly $\ln 2$. Our results provide a unified framework for generating and evaluating related binomial and Catalan-type series and reveal connections among central binomial coefficients, Catalan numbers, and harmonic numbers.
\end{abstract}



\bigskip

\newpage

\section{Introduction}

Infinite series for $1/\pi$ involving central binomial coefficients have a
	long history, extending well before Ramanujan's celebrated work on rapidly
	convergent series for $\pi$. Among the important early contributions are
	series due to A.\ R. Forsyth and J.\  W.\ L. Glaisher, now commonly referred to
	as Forsyth--Glaisher-type series. These identities typically
	involve powers of central binomial coefficients multiplied by rational
	functions of the summation index and provide a natural starting point for
	the study of many related series.
	
A classical example due to Forsyth is given in \cite{Forsyth}:
	\begin{equation}\label{Forsyth}
	\sum_{n=0}^{\infty}
	\frac{\binom{2n}{n}^{2}}{2^{4n}(2n-1)^2}
	=\frac{4}{\pi}.
	\end{equation}

Using elliptic functions, Glaisher proved in \cite{Glaisher} several
similar formulas, including the following identity:
	\begin{equation}\label{Glaisher}
	\sum_{n=0}^{\infty}
	\frac{\binom{2n}{n}^{2}}{2^{4n}(n+1)}
	=\frac{4}{\pi}.
	\end{equation}	
Glaisher subsequently obtained further series of this kind, including
identities with additional factors in the denominators. For historical
background and a discussion of the relationship between Forsyth's series
and Ramanujan's formulas, we refer the reader to \cite{Levrie-Nimbran}.

Further contributions to Forsyth--Glaisher-type series include the works of Almkvist \cite{Almkvist},  K\"u\c{c}\"uk and Sorgun \cite{Kucuk}, and Nimbran \cite{Nimbran}. Related developments in the broader
area of Ramanujan-type series for $1/\pi$ can be found in the works of
Bagis and Glasser \cite{Bagis}, 
Borwein et al. \cite{Borwein3},  Cooper \cite{Cooper}, Guillera 
\cite{Guillera4,Guillera5,Guillera6}, Guo and Zudilin
\cite{Guo}, Kouba \cite{Kouba}, Levrie \cite{Levrie}, Sun \cite{Sun}, Wan \cite{Wan}, Wang and Yang \cite{Wang},
Wei and Gong \cite{Wei}, and Zhang \cite{Zhang}.
	
Forsyth--Glaisher-type series should be distinguished from the classical
Ramanujan-type series for $1/\pi$. Although these two classes are related,
they have different historical origins and characteristic forms. In
particular, certain Forsyth--Glaisher-type series can be transformed into
Ramanujan-type series. This relationship motivates the investigation of
families of identities arising from the former class and the development
of elementary methods for deriving their generalizations.
	
In this paper, we develop an elementary approach to constructing and
evaluating parameterized families of Forsyth--Glaisher-type series for
$1/\pi$. Our approach is based on Newton's generalized binomial theorem,
identities for generalized binomial coefficients, and Euler's beta
integral. We first establish several auxiliary identities expressing
generalized binomial coefficients and products of odd integers in terms
of central binomial coefficients and Catalan numbers. These identities
allow us to transform suitable binomial series into infinite-series
evaluations involving squares of central binomial coefficients and
Catalan numbers.

We also establish additional families of series involving squares of
Catalan numbers. These identities arise from transformations connecting
generalized binomial coefficients, central binomial coefficients, and
Catalan numbers. In particular, Propositions~1--3 provide evaluations
involving products of odd factors, while their special cases yield
explicit formulas for series with squared Catalan numbers and
central binomial coefficients.

Another important part of the paper concerns series involving harmonic
and odd harmonic numbers. By differentiating suitable parameter-dependent
binomial identities, we derive several further families of infinite
series. These include formulas involving harmonic-number differences
and odd harmonic numbers, with explicit evaluations containing $1/\pi$
and logarithmic terms, particularly $\ln 2$. Related series involving
harmonic numbers and central binomial coefficients have been studied
by Campbell \cite{Campbell18,Campbell1,Campbell19},
Campbell et al. \cite{Campbell19-2}, Campbell and Cantarini
\cite{Campbell22}, Cantarini and D'Aurizio \cite{Cantarini}, Chen \cite{Chen}, Chu and Campbell \cite{Chu}, Hou et al.
\cite{Hou}, and Li and Chu \cite{Li}. Our results provide further
identities within this area through a systematic parameter-differentiation
approach.

We also consider series with shifted summation indices and additional
binomial factors. By combining finite-sum identities with transformations
of generalized binomial coefficients, we obtain further parameterized
formulas and their harmonic-number analogues. These results broaden the
range of series accessible through the proposed method.

The principal contribution of this work is a unified procedure for
constructing and evaluating several families of binomial and
Catalan-type series. The method brings together generalized binomial
identities, beta-integral evaluations, and parameter differentiation,
thereby revealing structural connections among central binomial
coefficients, Catalan numbers, and harmonic numbers.

The paper is organized as follows. Section~2 introduces the necessary
notation and preliminary identities, including several lemmas relating
generalized binomial coefficients to central binomial coefficients,
Catalan numbers, and products of odd integers. Section~3 develops the
main families of infinite-series identities, their special cases, and
their extensions involving harmonic and odd harmonic numbers.

\section{Preliminaries}

For nonnegative integers $i$ and $j$, the binomial coefficient is defined by
\begin{equation*}
\binom ij=
\begin{cases}
\dfrac{{i!}}{{j!(i - j)!}}, & \text{$i \ge j$};\\[2pt]
0, & \text{$i<j$}.
\end{cases}
\end{equation*}

For complex numbers $r$ and $s$, the generalized binomial coefficient is defined, wherever the expression is well defined, by
\begin{equation}\label{y89d722}
\binom rs= \frac{\Gamma (r + 1)}{\Gamma (s + 1)\Gamma (r - s + 1)},
\end{equation}
where $\Gamma$ denotes the gamma function. For $\Re(z)>0$, it is given by
$\Gamma(z)=\int_0^\infty e^{-t}t^{z-1}\,dt,
$ and it extends meromorphically to the complex plane, with poles at the nonpositive integers.

Relationships between binomial coefficients can be established by using \eqref{y89d722} together with the following identities:
\begin{equation}\label{Gamma12+}
\Gamma \Big( z + \frac{1}{2} \Big) = \frac{\sqrt \pi}{2^{2z}}\binom{{2z}}{z}\,\Gamma (z + 1),
\end{equation}
and, for every nonnegative integer $z$,
\begin{equation}\label{Gamma12-}
\Gamma \Big(- z + \frac{1}{2} \Big) = \frac{( - 1)^z 2^{2z}}{\binom{{2z}}{z}} \frac{\sqrt \pi}{\Gamma (z + 1)}.
\end{equation}
The first identity follows from the gamma-function duplication formula
\begin{equation}\label{dupl}
\Gamma(z+1)\Gamma\Big(z + \frac{1}{2} \Big) = \frac{\sqrt{\pi}}{2^{2k}}\,\Gamma (2z + 1),
\end{equation}
whereas the second follows from the reflection formula 
\begin{equation}\label{refl}
\Gamma\Big(\frac{1}{2}+z \Big)\Gamma\Big(\frac{1}{2}-z \Big) =(-1)^z\pi
\end{equation}
by taking $z$ to be a nonnegative integer and applying \eqref{Gamma12+}.

We also recall Euler's beta integral, which will be used later. For complex numbers $x$ and $y$ satisfying $\Re(x)>0$ and $\Re(y)>0$, we have
\begin{equation*}
\int_0^1 t^{x-1}(1-t)^{y-1}\,dt
=\frac{\Gamma(x)\Gamma(y)}{\Gamma(x+y)}.
\end{equation*}

We use the harmonic numbers $H_n$ and the odd harmonic numbers $O_n$, defined by
\begin{equation*}
H_n = \sum_{k=1}^{n} \frac{1}{k}, \qquad O_n = \sum_{k=1}^{n} \frac{1}{2k-1},
\end{equation*}
with $H_0=O_0=0$.

Another important ingredient in this paper is the sequence of Catalan numbers, which arises in numerous areas of mathematics, including lattice path enumeration, binary trees, polygon triangulations, noncrossing partitions, and generating function theory. The Catalan numbers are defined by
\begin{equation}\label{Cat-def-1}
C_n = \frac{1}{n+1}\binom{2n}{n},
\end{equation}
or, equivalently, by the recurrence
\begin{equation}\label{Cat-def-2}
C_0 = 1, \qquad C_n = \frac{2(2n-1)}{n+1}\,C_{n-1}, \qquad n\ge1.
\end{equation}

We conclude this section with several lemmas expressing generalized binomial coefficients and products of odd integers in terms of central binomial coefficients and Catalan numbers.
The preliminary identities established in this section serve two purposes. The first group expresses generalized
 binomial coefficients and products of odd integers in terms of Catalan numbers and central binomial coefficients.
The second group provides transformations that will be used to derive parameterized identities for series involving squared central binomial coefficients.
These relations form the algebraic basis of the series evaluations developed in Section~3.
\begin{lemma}\label{rkm6e32}
For nonnegative integers $r$ and $k$, we have
	\begin{equation}\label{Lemma1}
	\binom{r+\frac12}{k+1} = \frac{(-1)^{k+r}(2r+1)(r+1)!}{2^{2k+r+1}\prod\limits_{j=1}^{r}\bigl(2(k-j)+1\bigr)}\,C_kC_r.
	\end{equation}
Moreover, if $k\ge r\ge0$, then
	\begin{equation}\label{Lemma1_2}
	\binom{r+\frac12}{k+1} = (-1)^{k+r} \frac{2r+1}{2^{2k+1}} \frac{\binom{k}{r}}{\binom{2k}{2r}}\,C_k.
	\end{equation}
\end{lemma}
\begin{proof}
	By the definition of the generalized binomial coefficient,
	\[
	\binom{r+\frac12}{k+1}
	=\frac{1}{(k+1)!}
	\prod_{j=0}^{k}\Big(r-j+\frac12\Big).
	\]
	Separating the factors according to whether their indices are at most $r$ and simplifying the resulting products gives \eqref{Lemma1},
	where we have used \eqref{Cat-def-1}. 
	
For $k\ge r$, we also have
	\begin{equation}\label{prod-double}
	\prod_{j=1}^{r}\bigl(2(k-j)+1\bigr)
	=
	\frac{(2k-1)!!}{(2k-2r-1)!!}
	=
	\frac{(2k)!}{2^r(2k-2r)!}\frac{(k-r)!}{k!}.
	\end{equation}
	Substituting this expression and \eqref{Cat-def-1} into \eqref{Lemma1}, and simplifying, yields \eqref{Lemma1_2}.
\end{proof}
\begin{lemma}\label{vfidc73}
	If $m$ and $n$ are nonnegative integers with $m \le n$, then
	\begin{equation}
	\prod_{j = 1}^m \big(2(n - j) + 1\big) = \frac{m!}{2^m}\binom{n+1}{m}\frac{C_n}{C_{n - m}}.\label{lem2}
	\end{equation}
\end{lemma}
\begin{proof}
Using \eqref{prod-double}
and the definition \eqref{Cat-def-1} of the Catalan numbers, we obtain the stated identity after simplification.
\end{proof}
\begin{lemma}\label{zd6f3yu}
Let $r$ be a nonnegative integer and let $k\in\mathbb{C}$ be such that the generalized binomial coefficients in the following identity are well defined. Then
	\begin{equation}
	\binom{{k - \frac12}}{{k + r}} = \frac{( - 1)^r}{ 2^{ 2r+2k}} \frac{\binom{{2r}}{r} \binom{{2k}}{k}}{\binom{{k + r}}{r}}.\label{lem3}
	\end{equation}
\end{lemma}
\begin{proof}
	By \eqref{y89d722},
	\[ \binom{k-\frac12}{k+r} = \frac{\Gamma\!\left(k+\frac12\right)}{\Gamma(k+r+1)\Gamma\!\left(\frac12-r\right)}. 
	\]
	
	Using the reflection formula and \eqref{Gamma12+}, we obtain
	\[
	\binom{k-\frac12}{k+r}
	=
		\frac{(-1)^r\Gamma(2r+1)\Gamma(2k+1)}
	{2^{2r+2k}\Gamma(r+1)\Gamma(k+1)\Gamma(k+r+1)}.
	\]
	Finally, expressing the gamma functions in terms of generalized binomial coefficients gives \eqref{lem3}.
	\end{proof}
\begin{lemma}\label{x6lkv24}
	If $r$ and $k$ are nonnegative integers, then
	\begin{equation}\label{Lem4}
	\binom{{k}}{r + \frac12} = \frac{2^{r + 2k+ 2}}{(r + 2)! \,C_{r + 1}\binom{2k}{k}\pi}\prod_{j = 1}^r \big(2(k - j) + 1\big).
	\end{equation}
\end{lemma}
\begin{proof}
Using the gamma-function representation, we have
	\begin{equation*}
	\binom{k}{r+\frac12} = \frac{\Gamma(k+1)}{\Gamma\bigl(r+\frac32\bigl)\Gamma\big(k-r+\frac12\big)}.
	\end{equation*}
	Since
	\[
	\Gamma\Big(k+\frac12\Big) = \Gamma\Big(k-r+\frac12\Big) \prod_{j=1}^{r} \Big(k-j+\frac12\Big), 
	\]
	it follows that
	\begin{equation}
	\binom{k}{r+\frac12} = \frac{\Gamma(k+1)}{2^{r}\Gamma\big(k+\frac12\big)\Gamma\big(r+\frac32\big)}\prod_{j=1}^{r}\big(2(k-j)+1\big).
	\label{eq:proof1}
	\end{equation}

Applying the duplication formula and using
$$
\frac{\Gamma(k+1)}{\Gamma(k+\frac12)}
=
\frac{2^{2k}}{\sqrt{\pi}\binom{2k}{k}} \qquad\text{    and    }\qquad 
\Gamma\big(r+\frac32\big)
=
\frac{(2r+2)!\sqrt{\pi}}
{4^{r+1}(r+1)!},
$$
we obtain
\[
\binom{k}{r+\frac12}
=
\frac{2^{r+2}(r+1)!}
{\pi(2r+2)!}
\frac{2^{2k}}{\binom{2k}{k}}
\prod_{j=1}^{r}\bigl(2(k-j)+1\bigr).
\]
The definition of the Catalan numbers gives
$(r+2)!C_{r+1}
=\frac{(2r+2)!}{(r+1)!}$.
Thus the last expression simplifies to \eqref{Lem4}.
\end{proof}
\begin{lemma}\label{b5an9u4}
If $k$ and $r$ are nonnegative integers, then
\begin{equation*}
\binom{k + 1}{r + \frac32} = \frac{2^{r +2k+ 4}}{(r + 3)! \,C_{r + 2}C_k\,\pi}\prod_{j = 1}^r \big(2(k - j) + 1\big).
\end{equation*}
\end{lemma}
\begin{proof}
The result follows from Lemma~\ref{x6lkv24} by replacing $k$ with $k+1$ and $r$ with $r+1$, and then using \eqref{Cat-def-1}.
\end{proof}
\begin{lemma}\label{zgvhgwx}
If $r$ and $k$ are complex numbers such that the generalized binomial coefficients are defined, then
	\begin{equation}
	\binom{{k + r}}{{r + \frac12}} = \frac{{2^{2r +2k + 1} }}{{\pi ( 2r + 1)}}\frac{\binom{{k + r}}{r}}{\binom{{2r}}{r}\binom{{2k}}{k} }.\label{lem6}
	\end{equation}
\end{lemma}
\begin{proof}
	By \eqref{dupl},
	$$	\Gamma\Big(z+\frac12\Big)
	=\frac{\sqrt{\pi}}{2^{2z}}\binom{2z}{z}\Gamma(z+1).$$ Applying this identity with $z=r+1$ and $z=k$, and simplifying using
	$	\binom{2r+2}{r+1}
	=\frac{2(2r+1)}{r+1}\binom{2r}{r}$,
	we obtain \eqref{lem6}.
\end{proof}
\begin{lemma}\label{b819nks}
Let $k$ be a nonnegative integer and let $r\in\mathbb{C}$ be such that the generalized binomial coefficients in the following identity are defined. Then
	\begin{equation}\label{pz4a3ls}
	\binom{{r - \frac12}}{{k + r}} = \frac{( - 1)^k\binom{{2r}}{r}\binom{{2k}}{k}}{2^{2k+ 2r}\binom{{k + r}}{r}}.
	\end{equation}
\end{lemma}
\begin{proof}
The identity follows immediately from Lemma \ref{zd6f3yu} by interchanging $k$ and $r$.
\end{proof}
\begin{lemma}\label{vbyrvt6}
Let $k,r,m\in\mathbb{C}$ be such that the generalized binomial coefficients in the following identity are defined. Then
	\begin{equation}
	\binom{{k + r}}{{m + \frac12}} = \frac{1}{\pi }\,\frac{{2^{2m + 1 + 2k + 2r} }}{{(2m + 1)\binom{{2m}}{m}\binom{{2( {k + r})}}{{k + r}}}}\,\binom{{k + r - \frac12}}{m}.\label{lem8}
	\end{equation}
\end{lemma}
\begin{proof}
Put $z=k+r$. By Lemma~\ref{zgvhgwx},
\[
\binom{z}{m+\frac12}
=
\frac{2^{2z+1}}{\pi(2m+1)}
\frac{\binom{z}{m}}
{\binom{2m}{m}\binom{2(z-m)}{z-m}}.
\]
The gamma-function definition gives
$$\binom{z}{m}\binom{2(z-m)}{z-m}^{-1}
={2^{2m}}{\binom{2z}{z}}^{-1}
\binom{z-\frac12}{m}.$$
Substitution into the preceding identity yields \eqref{lem8}.
\end{proof}

\section{Families of Forsyth--Glaisher-type series}

In this section, we derive several families of Forsyth--Glaisher-type series by combining the auxiliary
identities established in Section~2 with an integral representation based on Newton's generalized binomial
theorem and Euler's beta integral. We first obtain evaluations involving squares of Catalan numbers and
products of odd integers. We then establish further identities involving squares of central binomial
coefficients and their parameterized extensions.
Finally, by differentiating suitable parameter-dependent
identities, we derive series involving harmonic and
odd harmonic numbers.
\begin{lemma}
	Let $p,r,s\in\mathbb{C}$ be such that the generalized binomial coefficients and all denominators below are well defined. Then
	\begin{equation}\label{g4how1g}
	\sum_{k = 0}^\infty  (-1)^k \frac{\binom{p}{k + 1}}{\binom{{k + r}}{s}} = \frac{1}{\binom{{r - 1}}{s}}  - \frac{s}{(r - s)\binom{{r + p - 1}}{{r - s}}}.
	\end{equation}
\end{lemma}
\begin{proof}
Multiplying the generalized binomial expansion
	\begin{equation*}
	\sum_{k = 0}^\infty  {( - 1)^k \binom{{p}}{{k + 1}}x^{k + 1} }  = 1 - \left( {1 - x} \right)^p 
	\end{equation*}
by $x^{r-s-1}(1-x)^{s-1}$ and integrating term by term over $[0,1]$, we obtain the result from Euler's beta integral.
\end{proof}
\begin{proposition} 
	If $m$ and $n$ are nonnegative integers, then
	\begin{align*}
	&\sum_{k = 0}^\infty  \frac{C_k^2}{2^{4k}  \prod\limits_{j = 1}^m \big(2(k - j) + 1\big) \prod\limits_{i = 1}^n \big(2(k - i) + 1\big)} \nonumber\\
	&\qquad\qquad =
	\frac{( - 1)^{m + n} 2^{m + n + 4}} {( m + 2)! (n + 2)! C_{m+1 } C_{n+1}}\left(\frac{2^{2m + 2n + 4}\binom{m+n+1}{m}}{\pi ( m + 1 ) (m+2)(n+2)\,C_{m+1}C_{n+1}}-1\right)\!.
	\end{align*}
\end{proposition}
\begin{proof}
Substituting $p=m+\frac12$, $r=1$, and $s=n+\frac32$ into \eqref{g4how1g} gives
	\begin{align*}
	\sum_{k = 0}^\infty  {( - 1)^k \frac{\binom{{m + \frac12}}{{k + 1}}}{\binom{{k + 1}}{{n +\frac 3 2}} }}  = ( - 1)^{n - 1} \frac{\pi }{2}(2n + 3) + \frac{{2n + 3}}{2n + 1}\frac{1}{\binom{{m + \frac12}}{{ m+n+1}}}.
	\end{align*}

Applying Lemmas \ref{rkm6e32} and \ref{b5an9u4}, together with the identity $$\prod_{j=1}^{m}\big(2(m+n-j)+1\big)=\frac{(m+n+1)! C_{m+n}}{2^m(n+1)! C_n},$$	
which follows from \eqref{lem2}, completes the proof.
\end{proof}

This proposition provides a two-parameter family of
series involving squares of Catalan numbers. Its
special cases yield explicit evaluations and illustrate
how products of odd integers can be incorporated into
Catalan-type series.
\begin{corollary}The following identities hold:	
	\begin{align}
	\sum_{k = 0}^\infty  \frac{C_k^2}{2^{4k}} &= \frac{{16}}{\pi } - 4\label{o33x7ff},\\
	\sum_{k = 0}^\infty  \frac{C_k^2}{2^{4k} \prod\limits_{j = 1}^n {\big(2(k - j) + 1\big)}}     &= 
	\frac{{( -2)^{n+ 3} }}{{(n + 2)!\,C_{n + 1} }}\left(1-\frac{{2^{2n + 2} }}{{ ( 2n + 1)C_n\pi }} \right)\!,\notag
	\end{align}
	and
\begin{align*}
		\sum_{k = 0}^\infty 
		\frac{C_k^2}{2^{4k}\prod\limits_{j = 1}^n \big(2(k - j) + 1\big)^2} 
		   = 
		\frac{2^{2n + 4}}{((n + 2)!)^2 C^2_{n+1}}\left(
		\frac{2^{4n + 3} }{(n + 2)(n+1)C_{n+1}\pi} - 1\right)\!.
		\end{align*}
\end{corollary}
\begin{proof}
		The first identity follows from the proposition with $m=n=0$.
	The second identity follows by setting $m=0$ and using \eqref{Cat-def-2}.
	
	Finally, setting $m=n$ in the proposition and expressing the resulting product in terms of Catalan numbers yields the third identity.
		\end{proof}

A useful feature of the parameter-dependent identities
	is that their differentiation produces additional factors involving harmonic numbers. This follows from
	the relation between derivatives of generalized binomial coefficients and the digamma function.
	For integer and half-integer arguments, these
	expressions can be rewritten in terms of ordinary harmonic numbers and odd harmonic numbers.
	We now apply this observation to derive several
	harmonic-number analogues of the preceding identities.
\begin{proposition}
Let $m$ and $n$ be nonnegative integers. Then
\begin{align*}
		\sum_{k = m}^\infty \frac{C_kC_{k-m}O_{k-m}}{2^{4k}\binom{k+1}{m}\prod\limits_{j = 1}^n{\!\big(2(k - j) + 1\big)}}  &= \frac{(-1)^{m-1}}{(m+1)(m+2)} \sum_{k = 0}^{m-1} \frac{(-1)^{k}\binom{m+2}{k+1}}{2^{4k} \prod\limits_{j = 1}^n \!\big(2(k - j)+1\big)}\frac{O_{m-k}C_k }{C_{m-k}} \\ 
		& + \frac{( - 1)^{n+m} 2^{2m + 3n + 7}\big( 4O_{m + 1}  - H_{m + n + 1}  - 2\ln 2\big)}{(m+1)( m + 2 )^2 n! (n+1)^2 ( n + 2)^2 }\frac{\binom{m+n+1}{n}}{C_{m + 1}^2 C^2_{n + 1}\pi }\\
		 & - \frac{{( - 1)^{n+m } 2^{ n + 4} }}{(m+1)(m+2)(n + 2)!}\frac{O_{m + 1}}{C_{m + 1} C_{n + 1}}.
		\end{align*}
\end{proposition}
	\begin{proof}
	Differentiate \eqref{g4how1g} with respect to $p$ and set $p=m+\frac12$, $r=1$, and $s=n+\frac32$ in~\eqref{g4how1g} to obtain, after some algebra,
	\begin{align*}
	\sum_{k = 0}^\infty  ( - 1)^k \frac{\binom{m + \frac12}{k + 1}}{\binom{k + 1}{n +\frac 32}} H_{m - k - \frac12}
	 = ( - 1)^{n - 1}  \frac{2n+3}{2}\pi H_{m + \frac12}  + \frac{{2n + 3}}{{2n + 1}}\, \frac{2H_{m + \frac12}  - H_{m + n + 1} }{\binom{{m +\frac 12}}{{m+n+1}}}.
	\end{align*}
	 Apply Lemmas~\ref{rkm6e32} and \ref{b5an9u4} and use the identity
$	 H_{r-1/2}=2O_{r}-2\ln 2$.
	\end{proof}
This result extends the preceding Catalan-number
	identities by incorporating odd harmonic numbers.
	It also illustrates how differentiation of the
	underlying parameter-dependent identity produces
	logarithmic contributions involving $\ln 2$.
\begin{corollary} The following evaluations  hold:
	\begin{gather*}
	\sum_{k = 0}^\infty  \frac{C_k^2  O_k}{2^{4k} }  = \frac{8}{\pi }\left( {3 - 2\ln 2} \right) - 4,\\
\sum_{k = 0}^\infty 
	\frac{C_k^2  O_k}{2^{4k}\prod\limits_{j = 1}^n (2(k - j) + 1)} 
	= \frac{( - 1)^n2^{n+3}}{(n+2)!C_n}\left( \frac{{2^{2n + 2} \left( {4 - H_{n + 1}  - 2\ln 2} \right)}}{{(n + 2)C_{n + 1}\,\pi  }} - 1\right)\!,
	\end{gather*}
and	
	\begin{align*}
		&(n^2+3n+2)^2 \sum_{k = 0}^\infty  \frac{C_k^2 O_{k}}{2^{4k}\binom{k+n+1}{n}^2}  \notag \\
		&\qquad\qquad=-2^{4n}\sum_{k = 0}^{n-1}  \frac{\binom{k+2}{n+1}^2}{2^{4k}} \frac{O_{n-k}}{C_{n-k}^2} +\frac{{2^{8n + 6} ( 4O_{n+1} - H_{2n + 1}  - 2\ln 2)}}{ (n+1)(n+2)C_{n + 1}^3 \pi} - \frac{2^{4n+4}\, O_{n+1}}{C_{n + 1}^2}.
		\end{align*}
\end{corollary}
\begin{proof}
	The first identity follows from the proposition by taking $m=0$.
	For the second identity, put $m=1$ and use 
	$$
	\prod_{i=1}^n \frac{1}{2(k-i)+1}=\frac{(-1)^{n+k}2^n\binom{n+2}{k+1}}{(n+2)!C_kC_{n-k}},\qquad k\leq n.
	$$
	
	Finally, setting $m=n$ in the proposition using 
	$$
	\prod_{i=1}^n \frac{1}{2(k+n-i)+1}=\frac{2^n (k+1)!C_k}{(k+n+1)!C_{k+n}}
	$$
	and simplifying gives the third identity.
\end{proof}
\begin{proposition}
	If $m$ and $n$ are nonnegative integers,  not both zero, then
\begin{align}\label{j66xa7q}
		&\sum_{k = 0}^\infty \frac{\binom{2k}{k}^2}{2^{4k}\prod\limits_{j = 1}^m \big(2(k - j) + 1\big) \prod\limits_{i = 1}^n \big(2(k - i) + 1\big)} 
		\notag\\
		&\qquad\qquad\qquad\qquad\qquad\qquad= \frac{(-2)^{3(n + m)}}{\pi m!n! (m+1)^2 (n+1)^2 (m+n)}\frac{\binom{m+n}{m}}{C^2_{n} C^2_{m}} .
		\end{align}
\end{proposition}
\begin{proof}
Setting $r=1$ in \eqref{g4how1g} and shifting the summation index, we obtain
	\begin{equation}\label{ngir6o2}
	\sum_{k = 0}^\infty  ( - 1)^k \frac{\binom{p}{k}}{\binom{{k}}{s}} = \frac{s}{(1 - s)\binom{p}{1 - s}}.
	\end{equation}
	Plugging $p=n-\frac12$ and $s=m+\frac12$ in \eqref{ngir6o2} gives
	\begin{equation}\label{pu8wcro}
	\sum_{k = 0}^\infty  ( - 1)^k \frac{\binom{n - \frac12}{k}}{\binom{k}{{m + \frac12}}}  =   \frac{2m+1}{2(m+n)}\frac{1}{\binom{{n - \frac12}}{{m + n}}} ,
	\end{equation}
	from which \eqref{j66xa7q} follows upon using Lemmas \ref{rkm6e32} and \ref{x6lkv24}.
\end{proof}

In contrast to Proposition~1, the present result is
	formulated directly in terms of squared central binomial
	coefficients. It provides a further family of
	Forsyth--Glaisher-type evaluations in terms of $1/\pi$.
\begin{corollary} We have
	\begin{equation*}
		\sum_{k = 0}^\infty \frac{\binom{{2k}}{k}^2}{2^{4k}(2k-1)}   = \frac{4}{\pi}
		\end{equation*}
	and
	\begin{equation}\label{series}
	\sum_{k = 0}^\infty  \frac{\binom{{2k}}{k}^2}{2^{4k}\prod\limits_{j = 1}^m (2(k - j) + 1)}
	  = \frac{{( - 2)^{3m} }}{{m!m\binom{{2m}}{m}^2 \pi }},\qquad m\geq1,
	\end{equation}
together with the special values
	\begin{gather}
	\sum_{k = 0}^\infty  
	\frac{\binom{{2k}}{k}^2}{2^{ 4k}(2k - 1) }
	  =  - \frac{2}{\pi }\label{qza7hvl},\\
	\sum_{k = 0}^\infty  \frac{\binom{{2k}}{k}^2}{2^{ 4k}(2k - 1)(2k - 3) }
	 = \frac{4}{{9\pi }},\label{qs50cr5}\\
	\label{pe3osl8}
	\sum_{k = 0}^\infty
	\frac{\binom{{2k}}{k}^2}{2^{ 4k}(2k - 3) } =  - \frac{{10}}{{9\pi }}.
	\end{gather}
\end{corollary}
\begin{proof}
	The first identity follows from \eqref{j66xa7q} by taking $m=n=1$, whereas the second one follows by setting $m=n$. The next two identities follow from \eqref{series} by taking $m=0$ and $m=1$, respectively.
		
Finally, \eqref{pe3osl8} follows from \eqref{qs50cr5} and \eqref{qza7hvl} by applying the partial-fraction decomposition
$
\frac{1}{(2k-1)(2k-3)}=
\frac{1}{2}\big(		\frac{1}{2k-3}-\frac{1}{2k-1}
		\big)$.
\end{proof}
\begin{proposition}
Let $m,n$ be nonnegative integers. Then
\begin{align}\label{Prop4new}
&\sum_{k = 0}^\infty  \frac{\binom{2k}{k}^2}{2^{ 4k}}\frac{O_{n}-O_{n-k}}{\prod\limits_{i = 1}^n\big(2(k - i) + 1\big)\prod\limits_{j = 1}^{m} \big(2(k - j)+1\big)} \nonumber\\
&\qquad\qquad\qquad\qquad\qquad = \frac{(-1)^{n+m} 2^{3(m + n)-1}}{ n!m!(n+m)}\frac{\binom{n+m}{n}}{\binom{2n}{n}^2 \binom{2m}{m}^2}\frac{H_{n+m-1}-2O_{n}+2\ln2}{\pi}.
\end{align}
\end{proposition}
\begin{proof}
Differentiating \eqref{ngir6o2} with respect to \(p\), we have
\begin{equation}\label{last}
			\sum_{k=0}^\infty
			(-1)^k \frac{\binom{p}{k}}{\binom{k}{s}}
			\left(H_p-H_{p-k}\right)
			=			
			\frac{s}{1-s}
			\frac{H_{p+s-1}-H_p}{\binom{p}{1-s}}.
	\end{equation}
			
Setting \(p=n-\frac12\) and \(s=m+\frac12\) in \eqref{last}, we get
				\begin{equation*}
\sum_{k=0}^\infty
			(-1)^k\frac{
			\binom{n-\frac12}{k}}{
			\binom{k}{m+\frac12}}
			\Big(H_{n-\frac12}-H_{n-\frac12-k}\Big)
			=		
			\frac{2m+1}{2m-1}
			\frac{H_{n-\frac12}-H_{n+m-1}}{\binom{n-\frac12}{\frac12-m}}.
			\end{equation*}
The desired result \eqref{Prop4new} now follows by applying
		Lemmas \ref{rkm6e32} and \ref{x6lkv24}, together with the
		identity
		$H_{n-\frac12}=2O_n-2\ln 2$.
	\end{proof}
	
This identity gives a harmonic-number analogue of the preceding central-binomial series and demonstrates how
parameter differentiation can be used to obtain additional evaluations involving harmonic-number differences.
\begin{corollary}
		We have
		\begin{gather}\label{Cor4new}
			\sum_{k = 0}^\infty  \frac{ \binom{2k}{k}^2}{{2^{  4k}(2k-1)}}  \,H_{\frac12-k}= \frac{8(\ln2-1)}{\pi},\\
		\sum_{k = 0}^\infty  \frac{ \binom{2k}{k}^2}{{2^{4k}(2k-1)^2}} \,H_{\frac12-k}= \frac{4(3-4\ln2)}{\pi},\nonumber\\
	\sum_{k = 1}^\infty   \frac{ \binom{2k}{k}^2}{{2^{  4k}(2k-1)}} \,{O_k} = \frac{2\ln2}{\pi},\nonumber
	\\
	\sum_{k = 1}^\infty 
	\frac{ \binom{2k}{k}^2}{{2^{  4k}(2k-1)(2k-3)}} \, O_{k} = -\frac{2(1+2\ln2)}{9\pi}.\notag
\end{gather}
\end{corollary}
\begin{proof}
The stated evaluations follow from \eqref{Prop4new} by setting $(m,n)=(0,1)$, $(1,1)$, $(1,0)$, and $(2,0)$, respectively, and using the evaluations
$$
		\sum\limits_{k = 1}^\infty 
		\frac{ \binom{2k}{k}^2}{{2^{4k}(2k-1)}} =1- \frac{2}{\pi},
		\qquad\quad
		\sum\limits_{k = 1}^\infty 
		\frac{ \binom{2k}{k}^2}{{2^{4k}(2k-1)^2}}
		= \frac{4}{\pi}-1,
		$$
	together with
		$H_{-k-\frac12}=2O_k-2\ln 2$, 
	and the special values
	$H_{\frac12}=2-2\ln 2$ 
	and $H_{-\frac12}=-2\ln 2.$
\end{proof}
\begin{theorem}\label{main_thm1}
If $m$ is a nonnegative integer and $r$ is a complex number that is not a negative integer and $r+m\ne0$, then
\begin{equation}\label{kgevius}
\sum_{k = 0}^\infty \frac{\binom{{2k}}{k}^2} {2^{4k}\binom{k + r}{r}\prod\limits_{j = 1}^m {{\big(2(k - j) + 1\big)}}}    = \frac{{( - 1)^m 2^{4r + 3m} }}{{m!( r + m)\binom{2( r + m )}{{r + m}}\binom{{2r}}{r}\binom{{2m}}{m}\pi }}.
\end{equation}
In particular,
\begin{equation}\label{u5u79zz}
\sum_{k = 0}^\infty \frac{\binom{{2k}}{k}^2} {2^{4k}\binom{{k + r}}{r}}     = \frac{{2^{4r} }}{{r\binom{{2r}}{r}^2 \pi }},\qquad 0\ne r\in\mathbb C\setminus\mathbb Z^{-},
\end{equation}
and
\begin{equation*}
\sum_{k = 0}^\infty  \frac{\binom{{2k}}{k}^2} {2^{4k}\prod\limits_{j = 1}^m\! \big(2(k - j) + 1\big)}   = \frac{{( - 8)^m  }}{{m!m\binom{{2m}}{m}^2 \pi }},\qquad m\geq1.
\end{equation*}
\end{theorem}
\begin{proof}
Identity \eqref{g4how1g} can be rearranged to read
\begin{equation}\label{flbygsu}
\sum_{k = 0}^\infty  {( - 1)^k \frac{{\binom{{m}}{k}}}{{\binom{{k + r}}{s}}}}  = \frac{s}{(r + 1 - s)\binom{{r + m}}{{r + 1 - s}}}.
\end{equation}
Setting $s=r+\frac12$ and writing $m-\frac12$ for $m$ in~\eqref{flbygsu} gives
\begin{equation*}
\sum_{k = 0}^\infty  {( - 1)^k \frac{{\binom{{m - \frac12}}{k}}}{{\binom{{k + r}}{{r + \frac12}}}}}  = \frac{2r + 1}{\binom{{r + m - \frac12}}{r + m - 1}},
\end{equation*}
from which \eqref{kgevius} follows by application of Lemmas \ref{rkm6e32} and \ref{zgvhgwx}.
\end{proof}

Theorem~5 provides a general parameterized extension of several classical Forsyth--Glaisher-type series.
	Its flexibility allows us to recover identities with
	different reciprocal products of odd factors and to
	derive further evaluations by suitable choices of
	the parameters. Moreover, differentiation with respect
	to the parameter $r$ leads to harmonic-number
	analogues, as demonstrated in Proposition~5.
\begin{remark}
	When $r$ is not a positive integer, \eqref{u5u79zz} does not give a Forsyth--Glaisher-like series. For example, for $r=\frac12$, we obtain the following classical identity:
	\begin{equation*}
	\sum_{k = 0}^\infty
	\frac{\binom{2k}{k}}{2^{2k}(2k+1)}
	= \frac{\pi}{2}.
	\end{equation*}
	This identity follows immediately from the Maclaurin expansion of the inverse sine function,
	\begin{equation*}
	\arcsin x =
	\sum_{k=0}^{\infty}
	\frac{\binom{2k}{k}}{2^{2k}(2k+1)}x^{2k+1},
	\qquad |x|\leq 1,
	\end{equation*}
	by setting $x=1$. 
\end{remark}
\begin{proposition}
If $r$ is a complex number that is not a non-positive integer, then
\begin{equation}\label{fv5v9y4}
\sum_{k = 0}^\infty 
\frac{\binom{{2k}}{k}^2}{2^{4k}\binom{{k + r}}{r}}
\,H_{k + r} = \frac{{2^{4r} }}{{r\binom{{2r}}{r}^2 \pi }}\left( {4\,O_r -H_{r-1}    - 4\ln 2} \right).
\end{equation}
In particular,
\begin{equation}\label{kvla2am}
\sum_{k = 0}^\infty \frac{\binom{2k}{k}^2}{2^{4k}(k + 1)}\,{H_{k}} = -\frac{16\ln2}\pi+4
\end{equation}
and
\begin{equation}\label{dynzpq2}
\sum_{k = 0}^\infty 
\frac{\binom{2k}{k}}{2^{4k}(2k + 1)}\,
O_{k + 1} = \pi \ln 2.
\end{equation}
\end{proposition}
\begin{proof}
Differentiate \eqref{u5u79zz} with respect to $r$. Identities \eqref{kvla2am} and \eqref{dynzpq2} are respective evaluations of \eqref{fv5v9y4} at $r=1$ and $r=\frac12$. Note that in obtaining the final form of \eqref{kvla2am}, we used \eqref{o33x7ff}.
\end{proof}
\begin{proposition}
If $m$ and $r$ are nonnegative integers, then
\begin{equation}\label{Prop6new}
\sum_{k = 1}^\infty 
\frac{\binom{2k}{k}^2}{\binom{k+r}{r}} \frac{ O_{m}-O_{m-k}}{2^{4k} \prod\limits_{j = 1}^m \!\big( 2(k - j) + 1\big)} = \frac{(-1)^m2^{4r + 3m-1}} {\pi m! (r +m) }\,\frac{H_{m+r-1}-2O_{m+r}+2\ln2}{\binom{2( r + m)}{r+ m}\binom{2r}{r}\binom{2m}{m}}.
\end{equation}
\end{proposition}
\begin{proof}
Differentiating \eqref{flbygsu} with respect to $m$, setting
$s=r+\frac12$, and then replacing $m$ by
$m-\frac12$, we obtain the desired result by applying Lemmas \ref{rkm6e32} and \ref{zgvhgwx}, together with the
identity
$H_{n-\frac12}=2O_n-2\ln 2$.
\end{proof}
\begin{corollary}The following series hold:
\begin{gather*}
\sum_{k = 1}^\infty 
\frac{ \binom{2k}{k}^2}{2^{4k}(2k-1)}\,  O_{k-1} = \frac{2(\ln2-1)}{\pi}+1,\\
\sum_{k = 1}^\infty 
\frac{ \binom{2k}{k}^2}{2^{4k}(k+1)}\,  O_{k} = \frac{4(1-\ln2)}{\pi},\\
\sum_{k = 1}^\infty 
\frac{ \binom{2k}{k}^2}{2^{4k}(k+1)(k+2)} \, O_{k} = \frac{8(5-6\ln2)}{27\pi},\\
\sum_{k = 1}^\infty 
\frac{ \binom{2k}{k}^2}{2^{4k}(k+1)(2k-1)}\,  O_{k-1} = \frac{4(6\ln2-11)}{9\pi}+1.
\end{gather*}
\end{corollary}
\begin{proof}
The results are obtained from \eqref{Prop6new} by setting $(r,m)=(0,1)$, $(1,0)$, $(2,0)$, and $(1,1)$, respectively, and using the evaluations
	\begin{gather*}
	\sum\limits_{k = 0}^\infty \frac{\binom{2k}{k}^2}{2^{4k}(2k-1)} = -\frac{2}{\pi},
	\qquad\qquad
	\sum\limits_{k = 0}^\infty 
	\frac{\binom{2k}{k}^2}{2^{4k}(k+1)}
	= \frac{4}{\pi},\\
	\sum\limits_{k = 0}^\infty 
	\frac{\binom{2k}{k}^2}{2^{4k}(k+1)(k+2)}
= \frac{16}{9\pi},\qquad\qquad
	\sum\limits_{k = 0}^\infty 
	\frac{\binom{2k}{k}^2}{2^{4k}(k+1)(2k-1)}
	 = -\frac{8}{3\pi},
	\end{gather*}
	together with
	$H_{-k-\frac12}=2O_k-2\ln 2
	$ and $H_{k-\frac32}=2O_{k-1}-2\ln 2
	$.
\end{proof}
\begin{theorem}\label{main_thm2}
If $m$ and $r$ are nonnegative integers such that $r>m$, then
\begin{align*}
\sum_{k = 0}^\infty  \frac{\binom{{2k}}{k}^2}{2^{4k}} \frac{{\prod\limits_{j = 1}^m {\left( {2(k + j) - 1} \right)} }}{{\prod\limits_{j = 1}^r {(k + j)} }}  = \frac{m!\,2^{4r - 3m} }{{r!( {r - m})}}\frac{\binom{2m}{m}}{\binom{{2( {r - m})}}{{r - m}}\binom{{2r}}{r}\pi}.
\end{align*}
In particular,
\begin{gather}
\sum_{k = 0}^\infty \frac{\binom{{2k}}{k}^2}{2^{4k}(k + 1)}   = \frac{4}{\pi }\label{f6865k6},\\
\sum_{k = 0}^\infty 
\frac{\binom{{2k}}{k}^2}{2^{4k}}
\frac{{2k + 1}}{{( {k + 1})( {k + 2})}}  = \frac{8}{{3\pi }}\label{j7t98o1},
\end{gather}
and 
\begin{equation}\label{mxsvkdf}
\sum_{k = 0}^\infty  
\frac{\binom{{2k}}{k}^2}{2^{4k}(k + 2)} 
 = \frac{{20}}{{9\pi }}.
\end{equation}
\end{theorem}
\begin{proof}
	Using the Pochhammer symbol $(a)_k$, we write the series as
	\begin{align*}
	\sum_{k=0}^{\infty}
	\frac{\binom{2k}{k}^2}{2^{4k}}
	\frac{\prod_{j=1}^m(2(k+j)-1)}
	{\prod_{j=1}^r(k+j)}
	=
	\frac{2^m(\frac12)_m}{r!}
	\sum_{k=0}^{\infty}
	\frac{(\frac12)_k(m+\frac12)_k}
	{(r+1)_k\,k!}.
	\end{align*}
Since $r>m$, Gauss's summation formula gives
	\[
	\sum_{k=0}^{\infty}
	\frac{(\frac12)_k(m+\frac12)_k}
	{(r+1)_k\,k!}
	=
	\frac{\Gamma(r+1)\Gamma(r-m)}
	{\Gamma(r+\frac12)\Gamma(r-m+\frac12)}.
	\]
	Hence, using the gamma-function duplication formula \eqref{dupl} and
	$(\frac12)_m=\frac{m!}{4^m}\binom{2m}{m}$,	we obtain
	\[
	\sum_{k=0}^{\infty}
	\frac{\binom{2k}{k}^2}{2^{4k}}
	\frac{\prod_{j=1}^m(2(k+j)-1)}
	{\prod_{j=1}^r(k+j)}
	=
	\frac{m!2^{4r-3m}\binom{2m}{m}}
	{r!(r-m)\binom{2(r-m)}{r-m}
		\binom{2r}{r}\pi},
	\]
	as claimed.
	
	Taking $(m,r)=(0,1)$ and $(1,2)$ gives
	\eqref{f6865k6} and \eqref{j7t98o1}, respectively. Finally,
	subtracting the identity for $(m,r)=(0,2)$ from
	\eqref{f6865k6}, and using
	$	\frac{1}{k+2}
	=\frac{1}{k+1}
	-\frac{1}{(k+1)(k+2)}$,
	yields \eqref{mxsvkdf}.
\end{proof}
\begin{lemma}\label{h1s6t7j}
If $n$ and $p$ are nonnegative integers and $m$ is a complex number such that $m+n\ne0$, then
\begin{equation}\label{sx3ccci}
\sum_{k = p}^\infty  {( - 1)^k \frac{{\binom{{n - \frac12}}{{k + n}}}}{\binom{k + n}{m + \frac12}}}  = ( - 1)^n \frac{{2m + 1}}{{2m + 2n}}\frac{1}{{\binom{{n - \frac12}}{{m + n}}}} - (-1)^n\sum_{k = 0}^{p + n - 1} {( - 1)^k \frac{\binom{{n - \frac12}}{k}}{\binom{k}{m + \frac12}}} .
\end{equation}
\end{lemma}
\begin{proof}
	Let
	$f(k)=(-1)^k{\binom{n-\frac12}{k}}{\binom{k}{m+\frac12}^{-1}}$.
	Applying the identity
	\[
	\sum_{k=p}^{\infty}f(k+n)
	=\sum_{k=0}^{\infty}f(k)-\sum_{k=0}^{p+n-1}f(k),
	\]
	valid for nonnegative integers $n$ and $p$ whenever the series converges, and then using \eqref{pu8wcro}, we obtain the  result.
\end{proof}
\begin{theorem}\label{ph75az1}
If $n$ and $p$ are nonnegative integers and $m$ is a complex number such that $m+n\ne0$ and $m\ne-\frac12$, then
\begin{align*}
&\sum_{k = p}^\infty  \frac{\binom{2k}{k}\binom{{2( k + n)}}{{k + n}}}{2^{4k}\binom{{k + n}}{n}\binom{{k + n - \frac12}}{m}}\nonumber\\
&\qquad\qquad   = \frac{ (-1)^n2^{4n+2m}}{\pi(m+n)\binom{2n}{n}\binom{2m}{m}}\left( \frac{{( - 1)^{m}2^{2(m+n)}\binom{{m + n}}{n} }}{{\binom{{2n}}{n} \binom{{2m}}{m} }} - \frac{2(m+n)}{2m+1}\sum_{k = 0}^{p + n - 1} {( - 1)^k \frac{{\binom{{n - \frac12}}{k}}}{{\binom{{k}}{{m + \frac12}}}}}\right)\!.
\end{align*}
\end{theorem}
\begin{proof}
The result follows by applying Lemmas \ref{b819nks} and \ref{vbyrvt6} to Lemma \ref{h1s6t7j}.
\end{proof}
\begin{corollary}
If $m$, $n$, and $p$ are nonnegative integers such that $m+n\ne0$ and $m\ne-\frac12$, then
\begin{align}\label{fds85y8}
&\sum_{k = p}^\infty  
\frac{\binom{2k}{k}\binom{2(k + n)}{k + n}}{\binom{k + n}{n}}\frac{2^{-4k}}{\prod\limits_{j = 1}^m \big( 2( k + n - j) + 1 \big)}
\nonumber\\
&\qquad  = \frac{(-1)^n2^{4n+m}}{\pi m!\binom{2n}{n}\binom{{2m}}{m}}\left(\frac{( - 1)^{m}2^{2n + 2m} }{m + n}\frac{{\binom{m + n}{n}}}{ {\binom{2n}{n} \binom{{2m}}{m}}} - \frac{2}{2m+1}\sum_{k = 0}^{p + n - 1} {( - 1)^k \frac{{\binom{n - \frac12}{k}}}{\binom{k}{m + \frac12}}}\right)\!.
\end{align}
\end{corollary}
\begin{remark}
Equation \eqref{series} is a special case of \eqref{fds85y8} being an evaluation at $n=p=0$.
\end{remark}
\begin{corollary}
If $n$ is a positive integer and $p$ is a nonnegative integer, then
\begin{equation*}
\sum_{k = p}^\infty   \frac{\binom{2k}{k}^2}{{2^{4k} \binom{n + k}{k}^2}}  = \frac{2^{4n-1}}{\pi \binom{{2n}}{n}^2}\left(\frac{{2^{4n } }}{n\binom{{2n}}{n}} -  \frac{4(-1)^n }{ 2n + 1}\sum_{k = 0}^{p + n - 1} {( - 1)^{k} \frac{{\binom{{n - \frac12}}{k}}}{{\binom{{k}}{{n + \frac12}}}}}\right)\!.
\end{equation*}
\end{corollary}
\begin{proof}
Set $m=n$ in~\eqref{sx3ccci} and use Lemmas \ref{zgvhgwx} and~\ref{b819nks}.
\end{proof}
\begin{theorem}\label{Th11}
	If $n$ and $p$ are nonnegative integers and $m$ is a complex number such that $m+n\ne0$ and $m\ne-\frac12$, then
\begin{align}
	&\sum_{k = p}^\infty  \frac{\binom{2k}{k}\binom{2( k + n)}{k + n}}{\binom{k + n}{n}}\frac{O_{m+1}-O_{k+n-m}}{2^{4k}\prod\limits_{j=1}^m \big(2(k+n-j)+1\big)}\nonumber\\
	&\qquad\quad  = \frac{(-1)^{m + n}2^{6n + 3m} }{\pi (m+n) m!}\,\frac{\binom{m+n}{n}} {\binom{2m}{m}^2 \binom{2n}{n}^2} \left(\frac{2n-1}{(2m+1)(m+n)} + H_{m+n} - 2O_m +2\ln2\right)\notag\\ 
	&\qquad\quad\quad- 2^{3n} n!\sum_{k = 0}^{p + n - 1}  \frac{\binom{2k}{k}^2}{2^{4k}} \frac{O_{m+1}-O_{k-m}}{\prod\limits_{i=1}^{m}\big(2(k-i)+1\big)\prod\limits_{j=1}^{n}\big(2(k-j)+1\big)}.\label{Eq_Th11}
	\end{align}
\end{theorem}
\begin{proof}
The result follows by  differentiating \eqref{sx3ccci} with respect to $m$ and simplifying.
\end{proof}
\begin{corollary} The following formulas hold:
\begin{gather*}
\sum_{k = 1}^\infty  \frac{\binom{2k}{k}^2}{ 2^{4k}(k+1)^2}\,O_{k}=\frac{24-16\ln2}{\pi}-4,\\
\sum_{k = 1}^\infty  \frac{\binom{2k}{k}^2}{2^{4k}(2k-1)}\,O_{k-1}=\frac{2\ln2-4}{\pi}+1,\\
\sum_{k = 1}^\infty  \frac{\binom{2k}{k}^2}{ 2^{4k}(2k-1)(2k-3)}\,O_{k-2}=\frac{26- 12\ln2}{27\pi}-\frac49.
\end{gather*}
\end{corollary}	
\begin{proof}
The stated identities follow from \eqref{Eq_Th11} by taking $(p,m,n)=(1,1,1)$, $(1,1,0)$, and $(1,2,0)$, respectively, and using 
	\begin{gather*}
	\sum\limits_{k = 1}^\infty \frac{\binom{2k}{k}^2}{2^{4k}(2k-1)} = 1-\frac{2}{\pi},
	\qquad
	\sum\limits_{k = 0}^\infty \frac{\binom{2k}{k}^2}{2^{4k}(k+1)^2} = -4+\frac{16}{\pi},\\
	\sum\limits_{k = 1}^\infty \frac{\binom{2k}{k}^2}{2^{4k}(2k-1)(2k-3)} =- \frac13+\frac{4}{9\pi},
	\end{gather*}
and applying  with
	$H_{-k-\frac12}=H_{k-\frac12}=2O_{k}-2\ln 2$.
\end{proof}

The results of this section establish several
	interconnected families of Forsyth--Glaisher-type
	series. Starting from parameterized binomial
	identities, we have obtained evaluations involving
	squared central binomial coefficients and Catalan
	numbers, as well as their harmonic and odd
	harmonic-number analogues. The various special cases
	demonstrate the flexibility of the method and
	illustrate how differentiation and parameter
	specialization lead to explicit evaluations
	involving $1/\pi$ and $\ln 2$.

\section{Conclusion}
In this paper, we have developed an elementary framework for constructing
and evaluating families of Forsyth--Glaisher-type series for $1/\pi$.
The method combines Newton's generalized binomial theorem, identities
for generalized binomial coefficients, Euler's beta integral, and
properties of the gamma and digamma functions. A key ingredient is the
conversion of generalized binomial coefficients and products of odd
integers into expressions involving central binomial coefficients and
Catalan numbers.

Using these tools, we have established several parameterized identities
for infinite series involving squared central binomial coefficients.
Their specializations yield numerous explicit evaluations in terms of
$1/\pi$, including formulas with reciprocal products of odd factors.
Further transformations lead to families of series involving squares
of Catalan numbers, providing alternative representations of the
resulting evaluations.
	
We have also shown that parameter differentiation provides a systematic
way to obtain harmonic-number analogues of the basic identities.
In particular, we derived series involving harmonic and odd harmonic
numbers, together with explicit evaluations containing logarithmic
terms, especially $\ln 2$. Additional identities with shifted summation
indices and finite-sum corrections further illustrate the flexibility
	of the approach.
	
The results demonstrate that a substantial collection of
Forsyth--Glaisher-type and related Ramanujan-like series can be derived
within a common framework based primarily on elementary binomial
identities and beta integrals. They also highlight the close
relationships among central binomial coefficients, Catalan numbers,
and harmonic numbers, showing how different types of series can arise
from the same underlying parameter-dependent identities.
	
The framework developed here suggests several directions for further
investigation. One possibility is to study higher derivatives of the
parameterized identities, which may lead to series involving generalized
harmonic numbers and higher-order harmonic sums. Another is to explore
other transformations of generalized binomial coefficients that could
produce new families of Catalan-type and binomial series, potentially
with evaluations involving additional logarithmic constants. These
questions may provide further applications of elementary methods to
the study of infinite series related to $1/\pi$.

\end{document}